
\input amstex
\documentstyle{amsppt}
\NoBlackBoxes

\define \U{U_q(\hat{\goth g})}
\define \a{\alpha}
\define \g{\gamma}
\define \ep{\epsilon}

\topmatter
\title On Drinfeld Realization of Quantum Affine Algebras\\
     \endtitle
\rightheadtext{Twisted quantum affine algebras}
\author Naihuan Jing
\endauthor
\address Department of Mathematics, North Carolina State University,
        Raleigh, NC 27695,  USA\endaddress
\abstract  
We provide a direct proof of the Drinfeld realization 
for the quantum affine algebras.
\endabstract

\thanks Supported in part by NSA grant MDA904-96-1-0087\endthanks
\endtopmatter

\document
\head 1. Introduction
\endhead

In 1987 Drinfeld \cite{Dr2} gave an extremely important realization of quantum
affine algebras \cite{Dr1}\cite{Jb}. This new realization has lead to numerous applications
such as the vertex representations \cite{FJ}\cite{J}. The proof of this
realization was not in print until Beck's braid group 
interpretation for the untwisted types \cite{B}. Some of lower rank
cases were also studied in \cite{D}\cite{S}. All these work started from the
quantum group towards the quantum loop realization, and were based
on Lusztig's theory of braid group action on the quantum enveloping
algebras \cite{L}. 
However, Drinfeld did give the exact isomorphism between two 
definitions of quantum affine algebras in \cite{Dr2}. In this paper we give
another proof directly from the Drinfeld isomorphism. Our proof
is self-contained and elementary and works from the opposite direction
from the quantum loop algebras towards the quantum groups.

In doing this, we discovered that there are rich structures held by the
$q$-loop algebra realization. We directly
deform the argument used by Kac \cite{K} to identify the affine Lie
algebras and the Kac-Moody algebra defined by generators and relations.
Here we must admit that the $q$-arguments are much more complicated than the 
classical analog, where we have the root space structure avaliable.  It is
nontrivial to properly deform usual brackets by $q$-brackets. 
As we have shown here in many cases there are strong indications 
for us to follow. The key are the following identities:
$$\align
[a, [b, c]_u]_v&=[[a,b]_x,c]_{uv/x}+x\, [b, [a,c]_{v/x}]_{v/x},
\qquad x\neq 0\\
[[a, b]_u, c]_v&=[a,[b, c]_x]_{uv/x}+x\, [[a, c]_{v/x}, b]_{u/x}, 
\qquad x\neq 0
\endalign
$$
where one needs to choose an appropriate $x$ to apply. 

We also discuss the Drinfeld realization for the twisted quantum 
affine algebras using the same approach. The results are used to construct the intertwining
operators between level one modules of twisted quantum affine algebras
in \cite{JK}.

I thank Professors Joseph Ferrar and Thomas Gregory for the organization of 
the conference on Lie algebras, where this paper was reported.

\head 2. Quantum affine algebras \endhead

Let $A=(a_{ij})$ $(i, j\in I=\{0, 1, \cdots, n\})$ be a generalized 
Cartan matrix of affine types \cite{K}.
Let ${\frak h}$ be a vector space over $\Bbb C(q^{1/2})$
with a basis $\{ h_0, h_1,
\cdots, h_n, d \}$ and define the linear functionals $\alpha_i \in
{\frak h}^*$ ($i\in I$) by
$$\alpha_i(h_j)=a_{ji}, \ \ \alpha_i(d)=\delta_{i,0} \ \ \textstyle {for} \
j\in I. \tag 2.1$$
Then the triple $({\frak h}, \Pi=\{\alpha_i|\ i\in I\},
\Pi^{\vee}=\{h_i|\ i\in I\})$ is the realization of the matrix $A$.
The Kac-Moody Lie algebra ${\frak g}$ associated with the matrix $A$ 
is called
the {\it affine Kac-Moody algebra of type $A$} (cf. [K]).
The algebra is generated as a Lie algebra by $e_i$, $f_i$, $h_i$ $(i\in I)$ and $d$ subject to the usual relations.
The elements of $\Pi$ (resp. $\Pi^{\vee}$) are
called the {\it simple roots} (resp. {\it simple coroots}) of ${\frak g}$.

The standard nondegenerate symmetric bilinear 
form $(\ |\ )$ on ${\frak h}^*$
satisfies
$$(\a_i|\a_i)=d_ia_{ij}, 
\ \ (\delta|\alpha_i)=(\delta|\delta)=0
\ \ \text {for all} \ i,j\in I. \tag 2.4$$
where 
$d_i=(\a_i|\a_i)/2$ are rational numbers given as follows:
$$\align
& (ADE)^{(1)}: d_i=1 \\
& C_n^{(1)}: d_0=1, d_i=1/2, d_n=1 \\
& B_n^{(1)}: d_0=d_i=1, d_n=1/2  \\
& G_2^{(1)}: d_0=d_1=1, d_2=1/3 \\
& F_4^{(2)}: d_0=d_1=d_2=1, d_3=d_4=1/2 
\endalign
$$
and $q_i=q^{d_i}$ for $i\in I$.
The {\it quantum affine Lie algebra $U_q(A)$}
is the associative algebra with 1 over $\Bbb C(q^{1/2})$
generated by the elements $e_i$, $f_i$ $(i\in I)$ and $q^h$ $(h\in P^{\vee})$
with the following defining relations :
$$
\aligned
&q^0=1, \ q^h q^{h'}=q^{h+h'} \ \ \text{for} \ h,h'\in P^{\vee},\\
&q^h e_{i} q^{-h}=q^{\alpha_i(h)} e_{i}, \ \
q^h f_{i} q^{-h}=q^{-\alpha_i(h)} f_{i} \ \ \textstyle {for} \ h\in P^{\vee}
(i\in I), \\
&e_{i}f_{j}-f_{j}e_{i} =\delta_{ij}
\dsize\frac {t_i-t_i^{-1}} {q_i-q_i^{-1}}, \ \textstyle {where} \
t_i=q^{h_i} \  \textstyle {and} \ i,j \in I, \\
&\sum_{m+k=1-a_{ij}} (-1)^m e_{i}^{(m)} e_{j} e_{i}^{(n)}=0,\\
&\sum_{m+n=1-a_{ij}} (-1)^m f_{i}^{(m)} f_{j} f_{i}^{(n)}=0
 \ \ \textstyle {for} \  i\neq j,
\endaligned \tag 2.5
$$
where $e_{i}^{(k)}=e_{i}^k/[k]_i !$, $f_{i}^{(k)}=f_{i}^k/[k]_i !$,
$[m]_i!=\prod_{k=1}^m [k]_i$, and
$[k]_i=\dsize\frac {q_i^k-q_i^{-k}} {q-q^{-1}}.$

Let $\Omega $ be the anti-algebra involution of $\U$ over $\Bbb C$ given by
$$
\Omega(e_i)=f_i, \quad \Omega(f_i)=e_i, \quad \Omega(q^h)=q^{-h}, 
\quad\Omega(q)=q^{-1}.
$$

Now we give the Drinfeld realization for the untwisted types.

Let {\bf U} be the associative algebra with 1 over $\Bbb C(q^{1/2})$
generated by the elements $x_i^{\pm}(k)$, $a_i(l)$, $K_i^{\pm 1}$,
$\gamma^{\pm 1/2}$, $q^{\pm d}$ $(i=1,2,\cdots,n, k\in \Bbb Z,
l\in \Bbb Z \setminus \{0\})$ with the following defining relations :
$$
\aligned
& [\gamma^{\pm 1/2}, u]=0 \ \ \text{for all} \ u\in \text{\bf U},\\
& K_iK_j=K_jK_i, \ \ K_iK_i^{-1} =K_i^{-1}K_i=1, \\
&[a_i(k), a_j(l)]=\delta_{k+l,0} \dsize\frac {[a_{ij}k]_i}{k}
\dsize\frac {\gamma^k-\gamma^{-k}}{q_j-q_j^{-1}},\\
&[a_i(k), K_j^{\pm 1}]=[q^{\pm d}, K_j^{\pm 1}]=0, \\
& q^d x_i^{\pm}(k) q^{-d}=q^k x_i^{\pm } (k), \ \
q^d a_i(l) q^{-d}=q^l a_i(l), \\
& K_i x_j^{\pm}(k) K_i^{-1}=q^{\pm (\alpha_i|\alpha_j)} x_j ^{\pm}(k), \\
& [a_i(k), x_j^{\pm} (l)]=\pm \dsize\frac {[a_{ij}k]_i}{k}
\gamma^{\mp |k|/2} x_j^{\pm}(k+l),\\
& x_i^{\pm}(k+1)x_j^{\pm}(l)-q^{\pm (\alpha_i|\alpha_j)} x_j^{\pm}(l)
x_i^{\pm}(k+1) \\
&= q^{\pm (\alpha_i|\alpha_j)} x_i^{\pm}(k) x_j^{\pm}(l+1)
-x_j^{\pm}(l+1) x_i^{\pm}(k), \\
& [x_i^{+}(k), x_j^{-}(l)]=\dsize\frac {\delta_{ij}}{q_i-q_i^{-1}}
\left( \gamma^{\frac {k-l}{2}} \psi_{i} (k+l)-\gamma^{\frac{l-k}{2}}
\varphi_{i} (k+l) \right), \\
& \text{ where $\psi_{i}(m)$ and $\varphi_{i}(-m)$ 
$(m\in \Bbb Z_{\ge 0})$
are defined by} \\
& \ \sum_{m=0}^{\infty} \psi_{i}(m) z^{-m}
=K_i \textstyle {exp} \left( (q_i-q_i^{-1}) \sum_{k=1}^{\infty} a_i(k) z^{-k}
\right),\\
& \ \sum_{m=0}^{\infty} \varphi_{i}(-m) z^{m}
=K_i^{-1} \textstyle {exp} \left(- (q_i-q_i^{-1}) \sum_{k=1}^{\infty} a_i(-k) z^{k}
\right),\\
&\ Sym_{l_1,\cdots, l_m}\sum_{s=1}^{m=1-a_{ij}}(-1)^s
\left[\matrix m\\s\endmatrix\right]_{i}
x^{\pm}_{i}(l_1)\cdots
x^{\pm}_i(l_s)\\
&\ \ \quad \cdot x^{\pm}_j(n)x^{\pm}_i(l_{s+1})\cdots
x^{\pm}_i(l_m)=0, \quad
\text{for $i\neq j$.}\\
\endaligned
\tag 2.11
$$

\proclaim {Lemma 2.1}
Let $I_0=\{1,2,\cdots,n\}$ be the index set for the simple
roots of a finite dimensional simple Lie algebra ${\frak g}_0$. 
Then for each $i\in I_0$, there exists a sequence of indices
$i=i_1, i_2, \cdots, i_{h-1}$ such that
$$
\aligned
& (\alpha_{i_1}|\alpha_{i_2})=\ep_{i_1}, \\
& (\alpha_{i_1}+\alpha_{i_2}| \alpha_{i_3})=\ep_{i_2}, \\
& \qquad\qquad\qquad\qquad \vdots \\
& (\alpha_{i_1}+ \cdots +\alpha_{i_{h-2}}|\alpha_{i_{h-1}})=\ep_{i_{h-2}},
\endaligned
\tag 2.12
$$
where $h$ is the Coxeter number of the Lie algebra ${\frak g}_0$,
and $\ep_i\in \Bbb Q_-$.
\hskip 1cm $\square$ \endproclaim

The twisted commutators $[b_1, \cdots, b_n]_{v_1\cdots v_{n-1}}$ and 
$[b_1, \cdots, b_n]'_{v_1\cdots v_{n-1}}$
is defined inductively by $[b_1, b_2]_v=[b_1, b_2]'_v=b_1b_2-vb_2b_1$ and 
$$\align
[b_1, \cdots, b_n]_{v_1\cdots v_{n-1}}&=[b_1, 
[b_2, \cdots, b_n]_{v_1\cdots v_{n-2}}]_{v_{n-1}}\\
[b_1, \cdots, b_n]'_{v_1\cdots v_{n-1}}&=[[b_1, 
\cdots, b_{n-1}]_{v_1\cdots v_{n-2}}, b_n]_{v_{n-1}}\endalign
$$
If $A$ is an antimorphism, then $A(
[b_1, \cdots, b_n]_{v_1\cdots v_{n-1}})=
[A(b_n), \cdots, A(b_1)]'_{v_{n-1}\cdots v_{1}}$. But if $B$ is an 
antimorphism such that $B(v_i)=v_i^{-1}$, then 
$$\align
& A([b_1, \cdots, b_n]_{v_1\cdots v_{n-1}})\\
&=[A(b_n), \cdots, A(b_1)]'_{v_{n-1}^{-1}\cdots v_{1}^{-1}}\\
&=v_1^{-1}\cdots v_{n-1}^{-1}[B(b_1), \cdots, B(b_n)]_{v_1\cdots v_{n-1}}
\endalign
$$

The following identities follow from the definition.
$$\align
[a, bc]_v&=[ab]_xc+x\, b[ac]_{v/x}, \qquad  x\neq 0\\
[ac, b]_v&=a[bc]_x+x\, [a, c]_{v/x}b, \qquad x\neq 0\\
\left[a, \left[b, c\right]_u\right]_v
&=\left[\left[a,b\right]_x,c\right]_{uv/x}
+x\, \left[b, \left[a,c\right]_{v/x}\right]_{v/x},
\qquad x\neq 0\tag2.13\\
\left[\left[a, b\right]_u, c\right]_v&=
\left[a,\left[b, c\right]_x\right]_{uv/x}
+x\, \left[\left[a, c\right]_{v/x}, b\right]_{u/x},
\qquad x\neq 0\tag2.14
\endalign
$$

In particular, we have
$$\gather
\left[a, \left[b_1, \cdots, b_n\right]_{v_1, \cdots, v_{n-1}}\right]=
\sum_{i}\left[b_1, \cdots, \left[a, b_i\right], \dots, b_n
\right]_{v_1, \cdots, v_{n-1}}. \\
[a, a, b]_{u\, v}=[a, a, b]_{v\, u}=a^2b-(u+v)aba+uv\,ba^2. \tag2.15
\endgather
$$
The Serre relation for the case of $A_{ij}=-1$ can be written as:
$$[x_i^{\pm}(m), x_i^{\pm}(m), x_j^{\pm}(n)]_{q_i, q_i^{-1}}=0,
\tag2.15'$$

\proclaim{Theorem 2.2} ({\rm [Dr]})
Fix an $\ep$-sequence $i_1, i_2, \cdots, i_{h-1}$, and let
$\theta=\sum_{j=1}^{h-1} \a_{i_j}$ be
the maximal root of the finite dimensional simple Lie algebra ${\frak g}$. 
Then there is a $\Bbb C(q^{1/s})$-algebra isomorphism
$\Psi: U_q({\frak g}) \to \textstyle {\bold U}$ defined by
$$
\aligned
& e_i \mapsto x_i^{+}(0), \ \ f_i \mapsto x_i^{-}(0), \ \
t_i \mapsto K_i \ \ \textstyle {for} \ i=1,\cdots, n, \\
& e_0 \mapsto [x_{i_{h-1}}^{-}(0), \cdots ,x_{i_2}^{-}(0),
x_{i_1}^{-}(1)]_{q^{\ep_1}\cdots q^{\ep_{h-2}}} \g K_{\theta}^{-1}, \\
& f_0 \mapsto a(-q)^{-\ep} \g^{-1}K_{\theta} 
[x_{i_{h-1}}^{+}(0), \cdots,
x_{i_2}^{+}(0), x_{i_1}^{+}(-1)]_{q^{\ep_1}\cdots q^{\ep_{h-2}}}\\
& t_0 \mapsto \gamma K_{\theta}^{-1}, \ \ q^d \mapsto q^d,
\endaligned
\tag 2.16
$$
where $K_{\theta}=K_{i_1}\cdots K_{i_{h-1}}$,
$h$ is the Coxeter number, $\ep=\sum_{i=1}^{h-2}\ep_i$,
and $s=1, 2, 3$, the quotient of long roots by
short roots. The constant $a$ is $1$
for simply types $A_n, D_n$, $a=[2]_1$
for $C_n^{(1)}$, and $a=[2]^{1-\delta_{1,i_1}}$ for $B_n^{(1)}$.
\endproclaim
\demo{Proof} Let $E_i, F_i, K_i, D$ be the images of $e_i, f_i, k_i, d$
in the algebra {\bf U}. We divide the proof into several steps.

Step 1. The elements $E_i, F_i, K_i, D$ satisfy the defining relations
of $\U$ given in (2.5). Clearly the defining relations of {\bf U}
imply that $E_i, F_i, K_i, i\neq 0$ generate a subalgebra isomorphic
to $U_q(\goth g)$. Thus we are left with relation involving $i=0$. For 
$i\neq 0$ we have
$$\align
\left[E_0, F_i\right]
&=\left[\left[x^-_{i_{h-1}}(0), \cdots, x_{i_1}^-(1)\right]_{q^{\ep_1}\cdots 
q^{\ep_{h-2}}}\g K_{\theta}^{-1}, x_i^-(0)\right]\\
&=-\left[x_i^-(0), x^-_{i_{h-1}}(0), \cdots, x_{i_1}^-(1)\right]_{q^{\ep_1}\cdots 
q^{\ep_{h-2}}q^{(\theta|\a_i)}}\g K_{\theta}^{-1}
\endalign
$$
We claim that 
$[x_i^-(0), x^-_{i_{h-1}}(0), \cdots, x_{i_1}^-(1)]_{q^{\ep_1}\cdots 
q^{\ep_{h-2}}q^{(\theta|\a_i)}}=0$ by the Serre relations. In fact this is
seen by looking at rank 3 cases. We show the argument by working out all 
cases of $A_3^{(1)}$. The first two use only one Serre relation, while the
third one uses two Serre relations.

$$\align
&[x^-_1(0), x^-_3(0), x_2^-(0), x_1^-(1)]_{q^{-1}\, q^{-1}\, q}\\
&=\left[x_3^-(0), \left[x_1^-(0), x_2^-(0), x_1^-(1)\right]_{q^{-1}\, q}
\right]_{q^{-1}}=0  
\endalign
$$
$$\align
&[x^-_3(0), x^-_3(0), x_2^-(0), x_1^-(1)]_{q^{-1}\, q^{-1}\, q}\\
&=\left[x_3^-(0), \left[\left[x^-_3(0), x_2^-(0)\right]_{q^{-1}}, x_1^-(1)
\right]_{q^{-1}}\right]_q\\  
&=\left[\left[x_3^-(0), x^-_3(0), x_2^-(0)\right]_{q^{-1}\, q}, x_1^-(1)
\right]_{q^{-1}}=0
\endalign
$$
$$
\align
&\left[x^-_2(0), x^-_3(0), x_2^-(0), x_1^-(1)\right]_{q^{-1}\, q^{-1}\, q^{-2}}\\
&=\left[x_2^-(0), \left[\left[x^-_3(0), x_2^-(0)\right]_{q^{-1}}, x_1^-(1)
\right]_{q^{-1}}\right]_{q^{-2}}\qquad\text{by (2.15')}\\  
&=\left[\left[x_2^-(0), x^-_3(0), x_2^-(0)\right]_{q^{-1}\, q^{-1}}, x_1^-(1)
\right]_{q^{-2}}\\
&\qquad
+q^{-1}\left[\left[x_3^-(0), x_2^-(0)\right]_{q^{-1}}, \left[x_2^-(0), x_1^-(0)
\right]_{q^{-1}}\right] \qquad\text{by (2.13)}\\
&=\left[\left[x_3^-(0), x_2^-(0), x_1^-(1)\right]_{q^{-1}\,q^{-1}}, x^-_2(0)
\right]_{q^{-2}} \qquad\text{by (2.14)}
\endalign
$$
which implies that $[x^-_2(0), x^-_3(0), x_2^-(0), 
x_1^-(1)]_{q^{-1}\, q^{-1}\, 1}=0$.

The Serre relations $\sum_{s=0}^{m=1-a_{ij}}(-1)^i
\left[\matrix m\\i\endmatrix\right]_ie_i^{s}e_je^{m-s}=0$
(resp. for $f_i$'s) for $i, j\neq 0$ are exactly the Serre relations 
in the Drinfeld realization. When $i$ or $j=0$ they boil down to the rank 
$3$ cases. We use $A_n^{(1)}$ to show the idea. 
$$\align
&e_0e_1^2-(q+q^{-1})e_1e_0e_1+e_1^2e_0\\
&=q^{-2}\left(\left[x_n^-(0), \cdots, x_1^-(1)\right]_{q^{-1}\cdots q^{-1}}x_1^+(0)^2\right.\\
&\qquad
-(q^2+1)x_1^+(0)\left[x_n^-(0), \cdots, x_1^-(1)\right]_{q^{-1}\cdots q^{-1}}\\
&\qquad \left.+q^2x_1^+(0)^2\left[x_n^-(0), \cdots, x_1^-(1)\right]_{q^{-1}
\cdots q^{-1}}\right)\gamma K_{\theta}^{-1}\\
&=q^{-2}
[x_1^+(0), x_1^+(0), x_n^-(0), \cdots, 
x_1^-(1)]_{q^{-1}\cdots q^{-1}\,1\,q^2} 
\qquad\text{by (2.15)}
\endalign
$$

Using commutation relations it follows that 
$$\align
&[x_1^+(0), x_1^+(0), x_n^-(0), \cdots, 
x_1^-(1)]_{q^{-1}\cdots q^{-1}\,1\,q^2}\\ 
&=\gamma^{-1/2}[x_1^+(0), x_n^-(0), \cdots, x_2^-(0),
K_1a_1(1)]_{q^{-1}\cdots q^{-1}\,q^2}\\ 
&=-[x_1^+(0), x_n^-(0), \cdots, x_3^-(0), 
x_2^-(1)K_1]_{q^{-1}\cdots q^{-1}\,q^2}\\ 
&=[x_1^+(0), x_n^-(0), \cdots, x_3^-(0),
x_2^-(1)]_{q^{-1}\cdots q^{-1}\,1}K_1=0 
\endalign
$$

Writing $\tilde{e_0}=e_0\gamma^{-1}K_{\theta}$, we have
$$\align
&e_1e_0^2-(q+q^{-1})e_0e_1e_0+e_0^2e_1\\
&=\left(x_1^+(0)\tilde{e_0}^2q-\left(1+q^{-2}\right)\tilde{e_0} x_1^+(0)\tilde{e_0}
+q^{-1}\tilde{e_0}^2x_1^+(0)\right)\gamma^2K_{\theta}^{-2}\\
&=q^{-1}[\tilde{e_0}, \tilde{e_0}, x_1^+(0)]_{1, q^2}\gamma^2K_{\theta}^{-2}\\
&=q^{-1}\left[\tilde{e_0}, \left[x_n^-(0), \cdots, x_2^-(0), -\gamma^{-1/2}
K_1a_1^+(1)\right]_{q^{-1}\cdots q^{-1}}\right]_{q^2}\gamma^2K_{\theta}^{-2}\\
&=q^{-1}\left[\left[x_n^-(0), \cdots, x_1^-(1)\right]_{q^{-1}\cdots q^{-1}},
\left[x_n^-(0), \cdots, x_2^-(1)\right]_{q^{-1}\cdots q^{-1}}\right]_{q}K_1
\gamma^2K_{\theta}^{-2}\\
&=-\left[\left[x_n^-(0), \cdots, x_2^-(1)\right]_{q^{-1}\cdots q^{-1}},
e_0\right]\gamma K_1K_{\theta}^{-1}\\
&=0
\endalign
$$
where we used our earlier result: $[e_0, f_i]=0$ for $i\neq 0$ and another
identity $[e_0, x_2(1)]=0$ by a similar argument as in (2.14).

Finally we check the relations $[e_i, f_i]=\frac{t_i-t_i^{-1}}{q_i-q_i^{-1}}
$. Again it suffices to see the case of $i=0$.
We want to give two cases to show the argument.

First we consider the case of $A_n^{(1)}$:
$$\align
&\left[\left[x^-_{n}(0), \cdots, x^-_{1}(1)\right]_{q^{-1} \cdots, q^{-1}},
\left[x^+_{n}, \cdots, x^+_{1}(-1)\right]_{q^{-1} \cdots, q^{-1}}\right]\\
&=\left[\left[\tilde{e_0}, x^+_{1}(0)\right], \cdots, x^+_{n}(-1)
\right]_{q^{-1} \cdots q^{-1}}\\
&\qquad\qquad +\left[x_n^+(0), \cdots, \left[\tilde{e_0}, x_1^+(-1)\right]
\right]_{q^{-1}\cdots q^{-1}}\\
&=\left[\left[x_{n-1}^-(0), \cdots, x_1^-(1)\right]_{q^{-1}\cdots q^{-1}}K_n,
x_{n-1}^+(0), \cdots, x_1^+(1)\right]_{q^{-1}\cdots q^{-1}}\\
&\qquad\qquad +\left[x_n^+(0), \cdots, x_2^+(0), K_1^{-1}\left[x_{n-1}^-(0), \cdots,
x_2^-(0)\right]_{q^{-1}\cdots q^{-1}}\right]\gamma\\
&=(-q^{-1})\left[e_0(n-1), f_0(n-1)\right]+q^{-1}K_1^{-1}\gamma\left[
x_{n-1}^-(0), \cdots, x_3^+(0), \right.\\
&\qquad\qquad \left.\left[x_2^+(0), x_{n-1}^-(0), \cdots, x_2^-(0)\right]_{q^{-1}\cdots q^{-1}\, 1}
\right]_{q^{-1}\cdots q^{-1}}\\
&=(-q)^{-n}\frac{\gamma K_{\theta}^{-1}K_n-\gamma^{-1}K_{\theta}}{q-q^{-1}}\\
&\qquad\qquad+(-q)^{-n+1}K_1^{-1}\cdots K_{n-1}^{-1}\gamma [x_n^+(0), x_n^-(0)]\\
&=(-q)^{-n}\frac{\gamma K_{\theta}^{-1}-\gamma^{-1}K_{\theta}}{q-q^{-1}}
\endalign
$$
where we have reasoned as follows:
the simplest Serre relations
$[e_0, e_i]=0$ for $i\geq 2$; an induction on rank $n$ as well as another
induction on the second bracket in line two. The elements $e_0(n-1), f_0(n-1)$
refer to the corresponding ones for $A_{n-1}^{(1)}$.

The computations in other cases are similar. I just give $C_2^{(1)}$ to
show some flavor. We write
$${\bar{e}}_0=[x_2^-(0), x^-_1(-1)]_{q^{-1}},
{\bar{f}}_0=[x_2^+(0), x^+_1(1)]_{q^{-1}} 
$$
$$\align
&\left[\left[x_1^-(0), x_2^-(0), x_1^-(1)\right]_{q^{-1}\, 1\, 1},
\left[x_1^+(0), x_2^+(0), x_1^+(-1)\right]_{q^{-1}\, 1\, 1}\right]\\
&=\left[\left[x_1^-(0), x_2^-(1)\right]_{q^{-1}}K_1, {\bar{f}}_0
\right]+\left[x_1^+(0), \left[x_1^+(0), \left[x_1^-(0),
 {\bar{e}}_0\right], {\bar{f_0}}\right]\right]\\
&=-\left[{\bar{e}}_0, {\bar{f_0}}\right]+
\left[x_1^+(0), \left[-x_2^+(-1)K_1^{-1},{\bar{e}}_0\right]\right]\\ 
&=q^{-1}\frac{\gamma K_1^{-1}K_2^{-1}-\gamma^{-1} K_1K_2}{q_1-q_1^{-1}}
+q^{-1}\left[x_1^-(0), x_1^+(0)\right]\gamma K_1^{-1}K_2^{-1}\\
&=q^{-1}[2]_1\frac{\gamma K_{\theta}^{-1}-\gamma^{-1} K_{\theta}}
{q-q^{-1}}.
\endalign
$$

Step 2. We show that the algebra $\bold U$ is generated by $E_i, F_i, 
K_i^{\pm 1}, D^{\pm 1}$. Write $\bold U'=
<E_i, F_i,K_i^{\pm 1}, D^{\pm 1}>$. 
The Cartan subalgebra is clearly generated by
$K_i^{\pm 1}$ and $D^{\pm 1}$. Rewriting (2.13) we have
$$\align
x_{i_1}^-(1)&= a[E_{i_2}, E_{i_3}, \cdots, E_{i_{h-1}}, E_0]_{q^{\ep_1}, 
\cdots, q^{\ep_{h-2}}} \\
x_{i_1}^+(-1)&= b[F_{i_2}, F_{i_3}, \cdots, F_{i_{h-1}}, F_0]_{q^{\ep_1}, 
\cdots, q^{\ep_{h-2}}} 
\endalign
$$
where $a, b$ are constants.
It then follows from Drinfeld relations that 
$$\align
a_{i_1}(1)&=K_{i_1}^{-1}\gamma^{1/2}[x^+_{i_1}(0), x_{i_1}^-(1)] \\
a_{i_1}(-1)&=K_{i_1}\gamma^{-1/2}[x^+_{i_1}(-1), x_{i_1}^-(0)]
\endalign
$$
which implies that $a_{i_1}(n)\in {\bold U}'$, subsequently 
$x_{i_1}^{\pm}(n)\in {\bold U}'$. Then we follow the $\ep$-sequence
and get that $x_j^{\pm}(n)\in {\bold U}'$. Thus ${\bold U}'={\bold U}$.

Step 3. We now have an algebra epimorphism $\Phi$: 
$U_q({\frak g}) \to \textstyle {\bold U}$. It is clear that $Ker\, \Phi=1$
when $q\to 1$ using Gabber-Kac \cite{GK}. Since quantization does not 
change the multiplicity we obtain that $\Phi$ is an automorphism.

\enddemo
$\blacklozenge$

\midinsert 
\topcaption{Table 2.1}
$\ep$-Sequences for simple Lie algebras
\endcaption

\centerline{\vbox{\tabskip=0pt\offinterlineskip
\def\tablerule{\noalign{\hrule}}
\halign to 360pt{\strut#& \vrule#\tabskip=0.1em plus0.1em&
      \hfil#& \vrule#& \hfil#\hfil& \vrule#&
      \hfil#& \vrule#\tabskip=0pt\cr\tablerule
&&\omit\hidewidth $\goth g$\hidewidth&&
  \omit\hidewidth $\ep$-Sequence. $\ep=\sum\ep_i$\hidewidth&&
  \omit\hidewidth $\ep$\hidewidth&\cr\tablerule
&& $A_n$ 
  && $\a_1 \overset{-1}\to{\rightarrow}\cdots \overset{-1}\to{\rightarrow} \a_n$
  && -n+1 &\cr\tablerule
&& $B_n$ 
  && $\a_1 \overset{-1}\to{\rightarrow}\cdots \overset{-1}\to{\rightarrow}
  \a_{n-1} \overset{0}\to{\rightarrow}\a_n \overset{-1}\to{\rightarrow}
  \cdots \overset{-1}\to{\rightarrow}\a_2$
 && -2n+4 &\cr\tablerule
&& $C_n$ 
  && $\a_1 \overset{-1/2}\to{\rightarrow}\cdots \overset{-1}\to{\rightarrow}
  \a_{n} \overset{-1/2}\to{\rightarrow}\cdots \overset{-1/2}\to{\rightarrow}
  \a_2 \overset{0}\to{\rightarrow}\a_1$
 && -n+1 &\cr\tablerule
&& $D_n$ 
  && $\a_1 \overset{-1}\to{\rightarrow}\cdots \overset{-1}\to{\rightarrow}
  \a_{n} \overset{-1}\to{\rightarrow}\a_{n-2} \overset{-1}\to{\rightarrow}
  \cdots \overset{-1}\to{\rightarrow}\a_2$
 && -2n+4 &\cr\tablerule
&& $E_6$ 
  && $\a_1 \overset{-1}\to{\rightarrow}\cdots \overset{-1}\to{\rightarrow}
  \a_{6} \overset{-1}\to{\rightarrow}\a_{3} \overset{-1}\to{\rightarrow}
  \a_2 \overset{-1}\to{\rightarrow}\a_4 \overset{-1}\to{\rightarrow}
  \a_3 \overset{-1}\to{\rightarrow}\a_6$
 && -10 &\cr\tablerule
&& $E_7$ 
  && $1\, 2\, 3\, 4\, 5\, 6\, 7\, 3\, 2\, 4\, 5\, 3\, 4\,
      7\, 3\, 2\, 1$, \quad $\ep_i=-1$
 && -16 &\cr\tablerule
&& $E_8$ 
  && $1\, 2\, 3\, 4\, 5\, 6\, 7\, 8\, 5\, 4\, 3\, 2\, 6\, 5\, 8\, 4\, 3\,
5\, 6\, 7\, 4\, 5\, 8\, 6\, 5\, 4\, 3\, 2\, 1$, \quad $\ep_i=-1$
 && -16 &\cr\tablerule
&& $F_4$ 
  && $1 \overset{-1}\to{\rightarrow}2 \overset{-1}\to{\rightarrow}3
  \overset{-1/2}\to{\rightarrow}4 \overset{-1/2}\to{\rightarrow}3 
  \overset{-1}\to{\rightarrow}2 \overset{-1}\to{\rightarrow}3  
  \overset{-1}\to{\rightarrow}4 \overset{0}\to{\rightarrow}3 
  \overset{-1}\to{\rightarrow}2 \overset{-1}\to{\rightarrow}1$
 && -7 &\cr\tablerule
&& $G_2$ 
  && $\a_1 \overset{-1}\to{\rightarrow}\a_2 
  \overset{-1/3}\to{\rightarrow}\a_{2} \overset{0}\to{\rightarrow}
  \a_1 \overset{-2/3}\to{\rightarrow}\a_2 $
 && -2 &\cr\tablerule
}}} 
\endinsert

\proclaim{Remark} If $\goth g$ is simply-laced, then $\ep=-h+2$. The
sequences are by no means unique, though $\ep$ is
independent from the choice of the sequences.
For example, we also have for $E_6$:
$$
\a_1 \overset{-1}\to{\rightarrow}\cdots \overset{-1}\to{\rightarrow}
\a_{6} \overset{-1}\to{\rightarrow}\a_{3} \overset{-1}\to{\rightarrow}
\a_2 \overset{-1}\to{\rightarrow}\a_4 \overset{0}\to{\rightarrow}\a_6  
\overset{-2}\to{\rightarrow}\a_3
$$
\endproclaim

\heading 3. Twisted Quantum Affine Algebras
\endheading

We now derive the Drinfeld realizations for the twisted types.

Let $X_N^{(1)}$ be a simply laced affine Cartan matrix. Let $\a_i'$ be the 
basis of the simple roots. The standard invariant bilinear form is
normalized 
as
$$(\a_i'|\a_i')=2r, \quad i=0, 1, \cdots, N$$

Let $L_q(X_N^{(1)})$ be the quantum affine algebra associated with $X_N^{(1)}$
realized in the Drinfeld quantum loop form. We denote the corresponding generators be putting an extra prime to distinguish.
Clearly the diagram automorphism $\sigma$ acts on the quantum affine algebra.
We will use a different indexing for the type $A_{2n}^{(2)}$ from \cite{K}.
The action of $\sigma$ is given as follows:
$$\align
& A_{N}: \sigma(i)=N-i\\
& D_{N}: \sigma(i)=i, 1\leq i\leq N-2; \sigma(N-1)=N\\
& E_6: \sigma(i)=6-i, 1\leq i\leq 5; \sigma(6)=6\\
& D_4: \sigma(1, 2, 3, 4)=(3, 2, 4, 1)
\endalign
$$
We construct a special invariant subalgebra ${\bold U}_{\sigma}$ generated
by:
$$\align
a_i(l)&=\frac{1}{[d_i]\sqrt r}\sum_{s=0}^{r-1}a_{\sigma^{s}(i)}(l)'\omega^{-ls},
\qquad K_i=\prod_{s=0}^{r-1}K_{\sigma^s(i)}'  \tag3.1\\
x_i^{\pm}(k)&=\frac{1}{[d_i]\sqrt r}\sum_{s=0}^{r-1}x_{\sigma^{s}(i)}(k)'
\omega^{-ks} \tag3.1 
\endalign
$$
where $\omega$ be a primitive $r$th root,
and
where $(d_0, \cdots, d_n)=(1, \cdots, 1, 2)$, $(1, 2, \cdots, 2, 1)$, 
$(2, 1, \cdots, 1, 1/2)$, $(1, 1, 1, 2, 2)$ and $(1, 1, 3)$,   
for $A_{2n-1}^{(2)}$, $D_{n+1}^{(2)}$, $A_{2n}^{(2)}$ , $E_6^{(2)}$
and $D_4^{(3)}$
respectively. We denote $[k]_j=
\frac{q^k_i-q^{-k}_i}{q_i-q^{-1}_i}$ if $j$ belongs to the $\sigma$-orbit of
$i$, then $[k]_j$ is defined for all $j=1, \cdots, N$ though we 
use only $[k]_i$ for $i\in \{0\}\cup\Gamma_{\sigma}=\{0, 1, \cdots, n\} $.
Easy and long calculation will lead to the following relations
presented in the Drinfeld realization of
twisted quantum affine algebras. 

\proclaim{Theorem 3.1}
The algebra 
${\bold U}_{\sigma}$ is the associative algebra with 1 over ${\Bbb C}(q^{1/2})$
generated by the elements $x_i^{\pm}(k)$, $a_i(l)$, $K_i^{\pm 1}$,
$\gamma^{\pm 1/2}$, $q^{\pm d}$ $(i=1,2,\cdots,N, k\in {\Bbb Z},
l\in {\Bbb Z} \setminus \{0\})$ with the following defining relations :
$$
\align
&x^{\pm}_{\sigma(i)}(k)=\omega^k x^{\pm}_i(k),  \ \ 
a_{\sigma(i)}(l)=\omega^l a_i(l) ,{}\\
&\  [\gamma^{\pm 1/2}, u]=0 \ \ \text{for all} \ u\in \text{\bf U},
{}\\
&\  [a_i(k), a_j(l)]=\delta_{k+l,0} 
\displaystyle\sum_{s=0}^{r-1}\frac {[k(\a'_i|\sigma^s(\a'_j))/rd_i]_i}{k}
\omega^{ks}\displaystyle\frac {\gamma^k-\gamma^{-k}}{q_j-q^{-1}_j},
\\
&\  [a_i(k), K_j^{\pm 1}]=[q^{\pm d}, K_j^{\pm 1}]=0, 
\\
&q^d x_i^{\pm}(k) q^{-d}=q^k x_i^{\pm } (k), \ \ 
q^d a_i(l) q^{-d}=q^l a_i(l), 
\\
&K_i x_j^{\pm}(k) K_i^{-1}=q^{\pm (\alpha_i|\alpha_j)} x_j ^{\pm}(k), 
\\
&\  [a_i(k), x_j^{\pm} (l)]=\pm \displaystyle\sum_{s=0}^{r-1}
\frac {[k(\a'_i|\sigma^s(\a'_j))/rd_i]_i}{k}\omega^{ks} 
\gamma^{\mp |k|/2} x_j^{\pm}(k+l),
\\
&\prod_s(z-\omega^sq^{\pm(\a'_i|\sigma^s(\a'_j))/r)}w)x_i^{\pm}(z)x^{\pm}_j(w)
=\prod_s(zq^{\pm(\a'_i|\sigma^s(\a'_j))/r}-\omega^s w)
x_j^{\pm}(w)x^{\pm}_i(z)
\\
&\  [x_i^{+}(k), x_j^{-}(l)]=\displaystyle
\sum_{s=0}^{r-1}\frac {\delta_{\sigma^s(i), j}\omega^{sl}}{q_i-q_i^{-1}}
\left( \gamma^{\frac {k-l}{2}} \psi_{i} (k+l)-\gamma^{\frac{l-k}{2}}
\varphi_{i} (k+l) \right), 
\endalign
$$
where $\psi_{i}(m)$ and $\varphi_{i}(-m)$ $(m\in {\Bbb Z}_{\ge 0})$
are defined by 
$$\align
&\sum_{m=0}^{\infty} \psi_{i}(m) z^{-m}
=K_i \textstyle {exp} \left( (q_i-q_i^{-1}) \sum_{k=1}^{\infty} a_i(k) z^{-k} 
\right),
\\
&\sum_{m=0}^{\infty} \varphi_{i}(-m) z^{m}
=K_i^{-1} \textstyle {exp} \left(- (q_i-q_i^{-1}) \sum_{k=1}^{\infty} a_i(-k) z^{k} 
\right),
\\
\text{Sym}_{z_1, z_2}&P_{ij}^{\pm}(z_1, z_2)\sum_{s=0}^{2}(-1)^s
\left[\matrix 2\\s\endmatrix\right]_{q^{d_{ij}}}
x^{\pm}_{i}(z_1)\cdots
x^{\pm}_i(z_s)x^{\pm}_j(w)x^{\pm}_i(z_{s+1})\cdot\\
\\ &\qquad\qquad\qquad\cdots x^{\pm}_i(z_2)=0, \quad
\text{for} A_{ij}=-1, \sigma(i)\neq j,\\
\ \textstyle{Sym}_{z_1, z_2, z_3}&\left[(q^{\mp 3r/4}z_1-{q^{r/4}+q^{-r/4})z_2+
q^{\pm 3r/4}z_3)x_i^{\pm}(z_1)}x_i^{\pm}(z_2)x_i^{\pm}(z_3)\right]=0,
\\
&\hskip 1in \text{for } A_{i,\sigma(i)}=-1\\
\endalign
$$
where Sym means the symmetrization over
  $ z_i $,
  $ P_{ij}^{\pm}(z,w) $
and
  $ d_{ij} $
are defined as follows:

$$\align
&\text{If }
  \sigma(i) = i, 
  \text{ then }
  P_{ij}^{\pm}(z,w) = 1 
  \text{ and }
  d_{ij}=r.
\\
& \text{If }
  A_{i,\sigma(i)} = 0
  \text{ and }
  \sigma(j) = j,
  \text{ then }
  P_{ij}^{\pm}(z,w) =
  \frac
     {z^r q^{\pm 2r}-w^r}
     { z q^{\pm 2} -w}
  \text{ and }
  d_{ij} = r.
\\
& \text{If }
  A_{i,\sigma(i)} = 0
  \text{ and }
  \sigma(j) \neq j,
  \text{ then }
  P_{ij}^{\pm}(z,w) = 1
  \text{ and }
  d_{ij} = 1/2.
\\
& \text{If }
  A_{i,\sigma(i)}= -1,
  \text{ then }
  P_{ij}^{\pm}(z,w) = 
  zq^{\pm r/2} + w 
  \text{ and } 
  d_{ij} = r/4. 
\endalign
$$ $\blacklozenge$
\endproclaim 

\midinsert 
\topcaption{Table 3.1}
$\ep$-Sequences for simple Lie algebras
\endcaption

\centerline{\vbox{\tabskip=0pt\offinterlineskip
\def\tablerule{\noalign{\hrule}}
\halign to 360pt{\strut#& \vrule#\tabskip=0em plus0.1em&
      \hfil#& \vrule#& \hfil#\hfil& \vrule#&
      \hfil#& \vrule#&
      \hfil#& \vrule#\tabskip=0pt\cr\tablerule
&&\omit\hidewidth $\goth g^{(r)}$\hidewidth&&
  \omit\hidewidth $\ep$-Sequence\hidewidth&&
  \omit\hidewidth $\ep$\hidewidth&&
   \omit\hidewidth $a$ \hidewidth&\cr\tablerule
&& $A_{2n-1}^{(2)}$ 
  && $\a_1 \overset{-1}\to{\rightarrow}\cdots \overset{-1}\to{\rightarrow}
  \a_{n-1} \overset{-2}\to{\rightarrow}\a_n \overset{-1}\to{\rightarrow}
  \cdots \overset{-1}\to{\rightarrow}\a_2$ && -2n+2 && $-2$ &\cr\tablerule
&& $D_{n+1}^{(2)}$ 
&& $\a_n \overset{-2}\to{\rightarrow}\cdots \overset{-2}\to{\rightarrow}\a_1$
 && -2n+2 && $(-2)^{n+1}$ &\cr\tablerule
&& $A_{2n}^{(2)}$ 
&& $\a_1 \overset{-1}\to{\rightarrow}\cdots \a_{n}\overset{0}\to{\rightarrow}
\a_{n-1} \overset{-1}\to{\rightarrow}\cdots \overset{-1}\to{\rightarrow}\a_2 
\overset{0}\to{\rightarrow}\a_1$
 && -2n+3 && $-[2]^{2n-2}$ &\cr\tablerule
&& $D_4^{(3)}$ 
  && $\a_1 \overset{-3}\to{\rightarrow}\a_2 \overset{-1}\to{\rightarrow}\a_1$
 && -4 && $3$ &\cr\tablerule
&& $E_6^{(2)}$ 
  && $\a_1 \overset{-1}\to{\rightarrow}\cdots \overset{-1}\to{\rightarrow}
  \a_{6}\overset{-1}\to{\rightarrow}\a_{3} \overset{-1}\to{\rightarrow}\a_2 
  \overset{-1}\to{\rightarrow}\a_4 \overset{-1}\to{\rightarrow}\a_3  
  \overset{-1}\to{\rightarrow}\a_6$
 && -10 && &\cr\tablerule
}}} 
\endinsert

\proclaim{Theorem 3.2} (\cite{Dr2})
Fix an $\ep$-sequence $i_1, i_2, \cdots, i_h$, and let
$\theta=\sum_{j=1}^{h-1} \a_{i_j}$ be
the maximal root of the finite dimensional simple Lie algebra ${\frak g}$. 
Then there is a $\Bbb C(q^{1/r})$-algebra isomorphism
$\Psi: U_q(\hat{{\frak g}}^{(r)}) \to \textstyle {\bold U}$ defined by
$$
\aligned
& e_i \mapsto x_i^{+}(0), \ \ f_i \mapsto \frac 1{p_i}x_i^{-}(0), \ \
t_i \mapsto K_i \ \ \textstyle {for} \ i=1,\cdots, n, \\
& e_0 \mapsto [x_{i_{h-1}}^{-}(0), \cdots ,x_{i_2}^{-}(0),
x_{i_1}^{-}(1)]_{q^{\ep_1}\cdots q^{\ep_{h-2}}} \g K_{\theta}^{-1}, \\
& f_0 \mapsto a(-q)^{-\ep} \g^{-1}K_{\theta} 
[x_{i_{h-1}}^{+}(0), \cdots,
x_{i_2}^{+}(0), x_{i_1}^{+}(-1)]_{q^{\ep_1}\cdots q^{\ep_{h-2}}}\\
& t_0 \mapsto \gamma K_{\theta}^{-1}, \ \ q^d \mapsto q^d,
\endaligned
\tag 2.13
$$
where $p_i=1$ for $\sigma(i)\neq i$, $p_i=i$ otherwise,
$K_{\theta}=K_{i_1}\cdots K_{i_{h-1}}$,
$h$ is the Coxeter number, $\ep=\sum_{i=1}^{h-2}\ep_i$,
and $d$ is a constant given by Table 3.1. 
\endproclaim

\demo{Proof} The proof is similar to that of the untwisted case.
\enddemo

\enddocument